\documentstyle[12pt]{article}

\newcommand{\be}{\begin{equation}}
\newcommand{\ee}{\end{equation}}
\newcommand{\zz} {\overline{z}}
\newcommand{\nn} {\noindent}
\begin{document}
\title{Bargmann representations for  deformed harmonic oscillators}
\author{Mich\`{e}le IRAC-ASTAUD and Guy RIDEAU\\ 
Laboratoire de Physique Th\'{e}orique de la mati\`{e}re condens\'ee\\
Universit\'{e} Paris VII\\2 place Jussieu F-75251 Paris Cedex 05, FRANCE}
\date{}
\maketitle

to be published in Reviews of Mathematical Physics

\begin{abstract}
Generalizing the case of the usual harmonic oscillator, we look for Bargmann representations corresponding to deformed harmonic oscillators. Deformed harmonic oscillator algebras are generated by  four operators $a, a^\dagger, N$ and the unity $1$ such as  $[a,N] = a ,\quad [a^\dagger ,N] = -a^\dagger$, $a^\dagger a = \psi(N)$ and $\quad aa^\dagger =\psi(N+1)$. We discuss the conditions of existence of a scalar product expressed with a true integral on the  space spanned by the eigenstates of $a$ (or $a^\dagger$). We give various examples, in particular we consider functions $\psi$ that are linear combinations of $q^N$, $q^{-N}$ and unity and that correspond to q-oscillators with Fock-representations or with non-Fock-representations. 
\end{abstract}

\section{Introduction}

\nn The harmonic oscillator Lie-algebra is defined by  four operators : the annihilation operator $a$, the creation operator $a^\dagger$ , the energy operator $N$ and the unity $1$  satisfying the following commutation relations :

\begin{equation}
[a,N] = a ,\quad [a^\dagger ,N] =
-a^\dagger 
\label{a1}
\end{equation}

\nn and

\be
[a,a^\dagger]=1
\label{usuel}
\ee
\nn where $a^\dagger$ is the adjoint of $a$ and $N$ is self-adjoint.

\nn This algebra has been deformed in many different ways (see in particular \cite{Bie}, \cite{Mac}, \cite{nous1}, \cite{nous2}, \cite{quesne}, \cite{Kos} ) and the representations of the deformed algebras were widely studied. In this paper, the deformed harmonic oscillator is defined by the  relations (\ref{a1}) and by the following relations between the three operators $a,a^\dagger$ and $N$:

\be
a^\dagger a = \psi(N) , \quad aa^\dagger =\psi(N+1)
\label{a3}
\ee

\nn where $\psi$ is a real analytical function.

\nn In the other formulations encountered in the literature, (\ref{a3}) is replaced by :

\be
[a,a^\dagger]=f(N,q)
\ee

\nn or

\be
[a,a^\dagger]_q\equiv aa^\dagger - q a^\dagger a=f_q(N,q)
\ee

\nn In these approach, the function $\psi$ is not given but results of solving the equations :

\be
f(N,q)= \psi (N+1) - \psi (N) \quad \mbox{or} \quad f_q(N,q)= \psi (N+1) - q \psi (N)
\ee

\nn The resolution of these equations leads to some arbitrariness that is eliminated in our formulation, $f$ and $f_q$ being uniquely determined in terms of $\psi$.

\nn Let us give some examples :

\nn - the usual harmonic oscillator defined by $f(N)=1$ corresponds to $
\psi_{usual}(N) = N +\sigma$.

\nn - the q-oscillator (\cite{Bie}, \cite{Mac}) defined by $f_q (N,q) = q^{-N}$ corresponds to 

\be
 \psi_{qosc} (N) = - q^{-N}/(q-q^{-1}) + \sigma q^N /(q-q^{-1}), \forall \sigma. 
\label{qosc}
\ee

\nn - the q-oscillator  defined by $f_q (N,q) = 1$ corresponds to 

\be
 \psi ^\prime_{qosc} (N) = (1-q)^{-1} + \sigma q^N , \forall \sigma.
\label{qosc2}
\ee

\nn - with the usual notation : 

\be
[x]=\frac{q^x -q^{-x}}{q-q^{-1}}
\ee

\nn  the function $\psi_{suq(2)} (N) = \sigma -[N-1/2]^2, \forall \sigma$,
 corresponds to $f(N) = -[2 N]$; that is to the deformation $su_q(2)$ of the Lie-algebra $su(2)$ after the identification $a=L_-, a^\dagger =L_+$ and $L_z= N$. 

\nn - $su_q(1,1) $ is obtained for the  $\psi_{suq(1,1)} (N)= - \psi_{suq(2)} (N)$.

\nn  Generalizing the pioneer work of Bargmann \cite{Bargmann} for the usual harmonic oscillator, the purpose of this paper is to study if  the deformed harmonic oscillator defined by (\ref{a1}) and (\ref{a3}) admits representations on one space of complex variable functions. In \cite{canada},  we restricted to the case where the function $\psi$ does not vanish. 
The scalar product of the representations we are looking for, is written with a true integral as in  (\cite{Bargmann}, \cite{canada}, \cite{jannussis}, \cite{azcarraga}, \cite{spiridonov}, \cite{kowalski} ) and contrarily to the works of (  \cite{gray}, \cite {chaichian}, \cite{bracken}, \cite{jurco}, \cite {quesne1}, \cite{odzij}, \cite{perelomov}) where a q-integration occurs. 

\nn In section 2, we recall how to build the irreducible representations on the basis of the eigenvectors of $N$. They are determined by the spectrum of $N$ which is depending on the  zeros of $\psi$. Then, in section 3, we discuss the existence of the coherent states that are defined as the eigenstates of the operators $a$ (or $a^\dagger$). 
In section 4, we study the possibility of  Bargmann representations. The formulation of the problem is done in a general framework. We show on various examples how the construction works  :  section 5 is devoted  to strictly positive function $\psi$ , other cases are considered in section 6.

\section{Representations}

\nn Let  $\mid 0>$ be  the  eigenvector of $N$ with eigenvalue $\mu$.
We built the 
 normalized vectors  $\mid n >$  

\be
\mid n> =\left\{
\begin{array}{ll}
\lambda _n a^{\dagger n}\mid 0>,& \quad n\in N^+\\
\lambda _n a^{-n}\mid 0>,& \quad n\in N^-
\end{array}
\right.
\ee

\nn with 

\be
\lambda_n^{-2} =\psi (\mu +n)!= \left\{\begin{array}{ll}
\prod_{i=1}^n \psi(\mu +i),& \quad n \in N^+\\
 \prod_{i=0}^{n+1} \psi(\mu +i),& \quad n \in N^-
\end{array}
\right.
\ee

\nn $N^+$ and $N^-$ are the set of integers $\geq 0$ and $<0$.

 \nn The vectors $\mid n>$ are the eigenvectors of $N$ with eigenvalue $\mu +n$ and span the Hilbert space $\cal{H}$. As $<n\mid aa^\dagger \mid n>$ is necessarily positive or zero, the construction of the increasing states  stops if it exists an integer $ \nu_+ $ such that 

\be
\psi(\mu+\nu_+ +1)= 0
\ee

\nn in which case the representation labelled by $\mu$ and $ \nu_+ $ admits a highest weight state $\mid \nu_+ >$. 
\nn We have an  analogous situation  for the decreasing states built with $a$,  when it  exists an integer $\nu_-$ such as  $\psi(\mu+ \nu_-)= 0$. The representation labelled by $\mu$ and $ \nu_- $ then admits a lowest weight state $\mid \nu_- >$.  

\nn We get different types of representations \cite{nous1},\cite{nous2},\cite {quesne}, \cite{Kos} :

\nn 1){\bf $\psi$ has no zero.}

\nn The inequivalent representations are labelled by the decimal part of $\mu$ and are defined by :

\be
\left\{
\begin{array} {ll}
a^\dagger \mid n >  =& (\psi (\mu +n+1))^{1/2} \mid n+1 > \\ 
a\mid n >  =& (\psi (\mu + n))^{1/2} \mid n-1 > ,\quad n \in Z \\ 
N \mid n >  =&  (\mu +n )\quad \mid n > 
 \end{array}
\right.
\ee

\nn The spectrum of $N$, $Sp N$, is $\mu + Z$. The operator $N$ has
 no lowest and no highest eigenstates. These representations, thus, are non equivalent to  Fock-representations and are called non-Fock-representations \cite{kulish} \cite {rideau}.  It is the case when $\psi$ is equal to $\psi_{qosc}$ with $q \in [0,1]$ and $\sigma \leq 0$. 

\nn An interpretation of this case \cite {kowalski} is obtained by identifying the states $\mid n >$ to the functions on a circle.

\nn {\bf 2) $\psi$ has zeros.}

\nn We are interested in the intervals where $\psi$ is positive :

\nn {\bf a) finite intervals }

\nn We can associate a representation to the intervals that have a length equal to an integer.

\nn The spectrum of $N$ is $[ \mu + \nu_-,\mu+\nu_+] \bigcap Z +\mu$.

\nn For example, in the case $\psi_{suq(2)}$, when $\sigma = [l+1/2]^2$, $l$ being a positive integer, the dimension of the representation is $2l+1$ and verifies

\be
\left\{
\begin{array}{ll}
a^\dagger\mid l,m>&= ([l+\frac{1}{2}]^2-[m+1-\frac{1}{2}]^2)^\frac{1}{2}\mid l,m+1>\\
a\mid l,m>&=([l+\frac{1}{2}]^2-[m-\frac{1}{2}]^2) ^\frac{1}{2}\mid l,m-1>\\
N\mid l,m>&=m\mid l,m>
\end{array}
\right.
\ee

\nn {\bf b) infinite intervals}

\nn The representations are similar to the  Fock-representation of the usual harmonic oscillator.

\nn The spectrum of $N$ is $ \mu + \nu_- + N^+$ or $\mu + \nu_+ + N^-$.

\nn Let us give an example :

\nn when $\psi$ is equal to $\psi_{qosc}$ with $\sigma = 1$, we recover the usual q-oscillator case such as

\be
\left\{
\begin{array} {ll}
a^\dagger \mid n >  =& [n+1]^{1/2} \mid n+1 > \\ 
a\mid n >  =& [ n]^{1/2} \mid n-1 > ,\quad n \in Z+  \\ 
N \mid n >  =&  n \quad \mid n >  
 \end{array}
\right.
\label{qosc1}
\ee

\nn The first step to build a Bargmann representation requires to study the coherent vectors.

\section{Coherent states}

 \nn We call coherent states \cite{Klauder}, the eigenvectors of the operator $a$ or $a^\dagger$. 

\nn The state $\mid z>= \sum_p c_p\mid p>$  is an eigenvector of $a$ if the coefficients $c_p$ verify the recursive relation 

\be
z c_p = \psi (\mu +p+1)^{1/2}c_{p+1}
\label{cp}
\ee

\nn - When the spectrum of $N$ is upper bounded, (\ref{cp}) implies that all the $c_p$ vanish and then that $a$ has no eigenvectors.

\nn If we look for the eigenvectors of $a^\dagger$, the situation is analogous : $a^\dagger$ has no eigenvectors if the spectrum of $N$ is lower bounded. 
Therefore, in the case (2.a) of the previous section as the spectrum of $N$ is finite, $a$ and $a^\dagger$ have no eigenvectors , hence  no Bargmann representation exists.

\nn - When the spectrum of $N$ is no upper bounded,
the  eigenvectors $\mid z>$ of $a$ take the form~:

\be
\left\{
\begin{array}{lll}
\mbox{when}& Sp N = Z+\mu,&\\
&&\\
\mid z>=&\sum_{n=-1}^{-\infty}z^{n}(\psi(\mu +n)!)^{1/2}\mid n>+\sum_{n=0}^{\infty}\frac{z^n}{ \psi(\mu +n)!^{1/2}}\mid n>,&\\
&&\\
\mbox{when}& Sp N = \nu_-+\mu +N^+,&\\
&&\\
\mid z>=&\sum_{n=-1}^{\nu_-}z^{n}(\psi(\mu +n)!)^{1/2}\mid n>+\sum_{n=0}^{\infty}\frac{z^n}{ \psi(\mu +n)!^{1/2}}\mid n>,& \nu_- < 0\\
&&\\
\mid z>=&\sum_{n=\nu_-}^{\infty}z^n (\psi(\mu +n)!)^{-1/2}\mid n>,&\nu_- \geq 0\\
\end{array}
\right.
\label{z}
\ee

\nn with the convention $\psi(\mu )! =1$.

\nn The domain $D$ of existence of the coherent states depends on the function $\psi$. Indeed, $\mid z>$ belongs to the Hilbert space spanned by the basis $\mid n>$ only if the series in the left hand side of (\ref{z}) are convergent in norm. 
\vspace{5mm}

\nn - When  $Sp N = Z+\mu$ , this implies that :

\be
\mid z \mid <lim_{p \rightarrow \infty} \psi(p)^{1/2}=r_2
\label{c1}
\ee

\nn and

\be
\mid z \mid >lim_{p \rightarrow -\infty} \psi(p)^{1/2}=r_1
\label{c2}
\ee

\nn Thus when $r_2=0$, the annihilation operator has no eigenvectors.

\nn When $r_1$ is smaller than $r_2$, the eigenstates of $a$ exist and their domain is $r_1<\mid z \mid < r_2$. When $r_1$ is larger than $r_2$, the annihilation operator $a$ has no eigenstates,  but then we can establish by analogous reasoning that the creation operator $a^\dagger$ has eigenstates if $r_1 \not= 0$.
\vspace{5mm} 

\nn - When the spectrum of $N$ is lower bounded, $Sp N = \mu + \nu_- + N^+$, the second condition (\ref{c2}) does not exist and  the eigenstates of $a$ always exist  provided $r_2 \not= 0$; their domain is defined by $\mid z \mid < r_2$.
\vspace{5mm}  

\nn When the spectrum of $N$ is upper bounded,  $Sp N = \mu + \nu_+ + N^-$, the eigenvectors of $a^\dagger$  exist only if $\mid z \mid < r_1$.

\nn To summarize, the eigenvectors of $a$ exist if :

\nn - $Sp N = \mu+Z$ and $r_2^2 \equiv \psi(+\infty)> r_1^2 \equiv \psi(-\infty)$, the domain of existence $D$ of the coherent states is $D = \{ z;r_1 < \mid z \mid < r_2 \}$ or 

\nn - $Sp N = \mu +\nu_- + N^+$, then $D = \{ z;\mid z \mid < r_2 \}$ . 

\nn The eigenvectors of $a^\dagger$ exist if :

\nn - $Sp N = \mu+Z$ and $r_2 < r_1$, then  $D = \{ z;r_2 < \mid z \mid < r_1 \}$ or

\nn - $Sp N = \mu +\nu_+ + N^-$, then $D = \{ z;\mid z \mid < r_1 \}$ .

\nn The part plays by $a$ and $a^\dagger$ being analogous, in the following we restrict to the case where the eigenstates of $a$ exist that is :

\nn -  $\psi $ is a  strictly positive function with $r_1 < r_2$

\nn - $\exists x_0$ such as $\psi (x_0)=0 $ and $\psi(x) > 0$ when $x > x_0$. 

\nn We do not study here the case where $r_1=r_2$.

\nn Although $\mu$ is a significant quantity as labelling inequivalent representations, it does not play a  part in the present problem. 
So we simplify the notation in assuming $\mu =0$ from now on. Indeed, this is equivalent to substitute $N-\mu$ to $N$ and $\psi_{\mu}(N)=\psi (\mu + N)$ to $\psi (N)$.

\section{ Bargmann representation}

\subsection{Representation space}

\nn Following the construction \cite{Bargmann}, in the Bargmann representation any state $\mid f>$ of $\cal{H}$ :

\be
\mid f>= \sum _{n\in Sp N}f_n \mid n> , \quad \sum _{n\in Sp N} \mid f_n \mid^2< \infty
\label{ff}
\ee

\nn is represented as the function of a complex variable $z$, $f(z) = < \zz \mid f >$, with a Laurent expansion   :

\be
f(z) = \left\{
\begin{array}{lll}
\sum _{n \geq 0} \frac{z^n f_n }{\psi (n)!^{1/2}} + \sum _{n < 0} z^n f_n (\psi (n)!)^{1/2},&
 Sp N = Z,&\\
&&\\
 \sum _{n \geq 0}\frac{z^n f_n }{\psi (n)!^{1/2}} + \sum _{n < 0}^{\nu_-} z^n f_n (\psi (n)!)^{1/2},
& Sp N = \nu_- + N^+,&\nu_- <0 ,\\
&&\\
\sum _{n \geq \nu_-} z^n f_n (\psi (n)!)^{-1/2},& Sp N = \nu_- + N^+,&
 \nu_- \geq 0 ,\\

\end{array}
\right.
\label{f}
\ee

\nn on the domain $D$ of definition of the eigenvectors of $a$, $r_1<\mid z \mid < r_2$. The space of the representation $\cal{S}$ is constituted with holomorphic functions in $D$, strongly depending on $\psi$ ( cf. subsection (4.3)). 

\nn In particular, to the basis vectors $\mid n> , n\in Sp N$, correspond the monomials~: 

\be
< \zz \mid n > = \left\{ \begin{array}{ll} 
z^n  (\psi (n)!)^{-1/2} , \quad &n \geq 0\\
z^n (\psi (n)!)^{1/2}, \quad &n <0
\end{array}
\right.
\label{n}
\ee

\subsection{Resolution of unity}

\nn A  Bargmann representation exists if  we can exhibit  a positive real function $F(x)$ such as 

\be
\int F(z \zz ) \mid \zz ><\zz \mid dz d\zz =1
\label{1}
\ee

\nn where the integration is extended to the whole complex plane and where $F(\mid z \mid ^2 )$ contains the characteristic function of the domain $D$ of existence of the  coherent states. Denoting $z= \rho e^{i\theta}$, Equation (\ref{1})reads :
\be
  \int_0^{2\pi}\frac{d\theta}{2}\int_{r_1^2}^{r_2^2} F(\rho ^2) d\rho^2 \mid \rho e^{-i\theta}>< \rho e^{i\theta} \mid =1
\label{01}
\ee

From the resolution of unity (\ref{1}), we obtain the scalar product in $\cal{S}$ :

\be
(g\mid f) = \int F(z \zz ) f(z)\overline{g(z)}dz d\zz
\label{prosca}
\ee

\nn and the representation of any linear operator $K$ of $ {\cal H}$ by  
 the kernel ${\cal K}(\zeta,\zz)= $ $<\overline \zeta \mid K \mid \zz >$ such as :

\be
(Kf)(\zeta) = \int F(z \zz ) {\cal K}(\zeta,\zz) f(z) dz d\zz 
\ee

\subsection{Reproducing Kernel}

\nn   Let $G(x)$ be the function  :

\be
G(x) =
\left\{
\begin{array}{lll}
 \sum _{n \geq 0} x^n  (\psi (n)!)^{-1} + \sum _{n < 0} x^n  \psi (n)!& Sp N=Z,&\\
\sum _{n \geq 0} x^n  (\psi (n)!)^{-1} + \sum _{n < 0}^{\nu_-} x^n  \psi (n)!,& Sp N=\nu_- +N^+,&\nu_-<0\\
\sum _{n \geq \nu_-} x^n  (\psi (n)!)^{-1},& Sp N=\nu_- +N^+,&\nu_-\geq 0\\
\end{array}
\right. 
\label{G}
\ee

\nn We easily verify that $G(\zeta z) = <\zz \mid \zeta>$, so that $G(\zeta z) $ is the function of $\cal{S}$ representing the coherent state $\mid \zeta>$. The function $G(x)$ plays a prominent part for, if (\ref{1}) should be true, we will have
 
\be
\int F(z \zz ) < \overline{\zeta}\mid \zz ><\zz \mid f > dz d\zz = < \overline{\zeta}\mid f >
\ee

\nn and  $ G(\zeta \zz)$ appears to be  one reproducing kernel.

\nn Moreover, we get from the Schwarz inequality :

\be
\mid f(z) \mid \leq \mid f\mid ^{1/2} \mid G(z \overline{z})\mid ^{1/2}
\label{<}
\ee

\nn Thus,   $\cal{S}$ is the set of holomorphic functions the growth of which on the edge of  $D$, is controlled by the growth of $\mid G(x) \mid^{1/2}$.

 \nn We easily prove  that, in  the Bargmann representation, $a^\dagger$  is the multiplication by $z$ and $N$ is the operator $zd/dz$, as in the usual case, while  $a$ is the operator $z^{-1}\psi(zd/dz)$. As $G(\zeta z)$  corresponds  to the coherent state $\mid \zeta >$, we obtain for $G$ the following equation
\be
x G(x) = \psi(x \frac{d}{dx}) G(x)
\label{GG}
\ee

\nn which could be obtained directly from the expansion (\ref{G}).

\subsection{ Weight Function}

\nn Let us introduce the Mellin transform $\hat{F}(\rho)$ of the weight function $F(x)$ :

\be
\hat {F}(\rho)=\int_0^{\infty} F(x) x^{\rho -1} dx = \int_{r_1^2}^{r_2^2} F(x) x^{\rho -1} dx
\label{mel}
\ee

\nn From (\ref{n}) and (\ref{1}), we deduce that $F(x)$ must verify the following condition~:

\be
\hat{F}(n+1)= \left\{ \begin{array}{lll}
 \psi(n)!, & n\geq 0&\\
 (\psi(n)!)^{-1}, & n <0& , n \in Sp N
\end{array}
\right.
\label{moments}
\ee

\nn Let us remark that 

\be
\hat {F}(\rho) \leq \hat {F}(n)+\hat {F}(n+1), \quad n \leq Re \rho <n+1
\label{nn}
\ee

\nn as $F(x)$ is a positive function.
As the Mellin transform of $F$ exists for all the integers belonging to the spectrum of $N$, then, due to (\ref{nn}), it exists for any $\rho$  such as $Re \rho$ is greater than the lowest bound $\nu_-$ of $Sp N$, $\nu_-$ can be finite or $-\infty$. 

\nn Formula (\ref{moments}) is equivalent to

\be
\hat{F}(n+1)=\psi (n) \hat{F}(n),\mbox{with} \quad \hat{F}(1)=1
\label{M}
\ee

\nn which ensures that the operators $a^\dagger = z$,  $a = z^{-1}\psi(zd/dz)$ be adjoint on the basis $\mid n>$.

\nn The function $\psi$ being given,  $\hat{F}$ verifying (\ref{M}) must be interpolated in order to get $F$ as the Mellin inverse of (\ref{mel}). 

\nn Equation (\ref{M}) can be obviously interpolated by :

\be
\hat {F}(\rho+1)= \psi(\rho)\hat {F}(\rho)
\label{mel1}
\ee

\nn Any other interpolation  of $\hat{F}(n)$ is obtained by :

\be
\hat {F}_{h}(\rho)= \hat {F}(\rho)\hat{h}(\rho)
\ee

\nn where $\hat {F}(\rho)$ is some solution of (\ref{mel1}) and where $\hat{h}(\rho)$ is equal to $1$ on the spectrum of $N$.

\nn In order to tackle the discussion on the existence of a Bargmann representation, we first state the following lemmas :

\nn \underline{Lemma 1} : If there exists one particular solution of (\ref{mel1}) which does not admit a Mellin inverse, the same holds for every solution of (\ref{mel1}).

\nn \underline{Proof} : Two solutions of (\ref{mel1}) differ by a multiplicative factor $\hat{k}(\rho)$ which is a periodic function of period $1$ :

\be
\hat{k}(\rho +1)=\hat{k}(\rho)
\label{k}
\ee

\nn Let $\hat {F}(\rho)$ be one particular solution of (\ref{mel1}), without Mellin inverse. We cannot find $\hat{k}(\rho)$ such as $\hat{k}(\rho)\hat {F}(\rho)$ has a Mellin inverse. Indeed, the behavior of $\hat{k}(\rho)$ at infinities on a parallel to the imaginary axis must compensate the corresponding bad behavior of $\hat {F}(\rho)$. This implies that $\hat{k}(\rho)$ has a Mellin inverse $k(x)$. Then from Equation (\ref{k}), we obtain that $k(x)$ must verify

\be
x k(x) = k(x)
\ee

\nn and then $k(x)$ is equal to $ k \delta (x-1)$, so that $\hat{k}(\rho)$ is a constant in contradiction to the assumption that $\hat{k}(\rho)\hat {F}(\rho)$ have a Mellin inverse.

\nn \underline{Lemma 2} : If $\hat{h}(\rho)$ is equal to $1$ on $Sp N$ and admits a Mellin inverse $h(x)$, $\hat{h}(\rho)$ is the identity.

\nn \underline{Proof} : Let us assume that $Sp N$ contains $N^+$. Then we have :

\be
\hat{h}(n+1)= \int_0^\infty h(x) x^n dx =1
\ee

\nn From this equation, we easily deduce :

\be
\int_0^\infty h(x) e^{\rho x} dx = e^\rho
\ee

\nn so that $h(x) =\delta (x-1)$ and $\hat{h}(\rho)=1$.

\nn The same method works when $Sp N = \nu_- + N^+, \nu_->0$.

\nn \underline{Discussion} :

\nn 1) If there exists a solution of (\ref{mel1}) with a positive Mellin inverse, our problem is solved.

\nn 2) If we exhibit a solution of (\ref{mel1}) without Mellin inverse, we know from Lemma 1 that no solution of (\ref{mel1}) with Mellin inverse exists and from Lemma 2 that these solutions cannot be corrected by a suitable multiplicative factor. 

\nn 3) If there exists a solution of (\ref{mel1}) with a Mellin inverse not strictly positive on $D$, we cannot conclude, for we cannot ensure this solution cannot be improved by a suitable multiplicative factor leading to positivity.

\nn Therefore, a Bargmann representation can be obtained only if there exist solution of (\ref{mel1}) with Mellin inverse, so that the interpolation problem (\ref{M}) is completely solved by (\ref{mel1}) precisely.

\nn Concerning the remaining positivity condition, let us point out that it will be satisfied if $\psi(\rho)$ is such that $\psi(\rho +i _sigma)$ is a definite positive function of $\sigma$ for any $\rho \in Sp N$. This can be proved by an adaptation to the Mellin transform of the well-known Bochner theorem for Fourier transform. But this condition on $\psi(\rho)$ is only a sufficient condition as it is readily seen from the counter-example of the usual oscillator : $\psi(\rho)= \rho$.

\nn Therefore, the only interpolation of (\ref{M}) to be considered in the following is precisely the simplest one, namely (\ref{mel1}).

\nn Remains one case where we cannot conclude, namely when the Mellin inverse of  the solution of (\ref{mel1} exists but is not positive on D, we cannot ensure that the result can be improved by a suitable factor leading to the positivity.

\nn To summarize, a Bargmann representation can be defined on a deformed harmonic oscillator algebra  if the coherent states exist and if it exists at least one solution of (\ref{mel1}) with a positive Mellin inverse on D.

\nn Furthermore, let us remark that (\ref{mel1}) can be written :

\be
\int F(x) x^{\rho} dx = \int F(x) \psi(x\partial _x +1) x^{\rho -1} dx
\ee

\nn The right member is equal to $\int \left( \psi(-x\partial_x) F(x)\right) x^{\rho -1} dx$ if the integration by parts (or the change of variables) can be done without extra terms and this gives :

\be
x F(x) = \psi(-x \frac{d}{dx}) F(x)
\label{FF}
\ee

\nn This equation, when it holds, can be used to study the positivity of $F(x)$, when we cannot obtain an explicit expression for this function. 

\nn In the next section \cite{canada} , we restrict to strictly positive  $\psi$ and to $F$  solution of (\ref{FF}), fast decreasing at infinity and at the origin. 

\nn  In section 6, we illustrate the construction of Bargmann representation  with specific vanishing functions $\psi$, in particular we discuss the cases of the q-oscillators, looking for a resolution of the identity involving a true integral.

\section {$\psi$ strictly positive }

 \nn In this section, the eigenvalues of $N$ are all the integers. The coherent states lie in the whole complex plane, except the origin, in the two first subsections and in a ring in the last one. 

\subsection{ the q-oscillator}

\nn The function $\psi$ associated to the q-oscillator (\ref{qosc}) is strictly positive if $q<1$ and $ \sigma \leq 0$. We easily see that when $\sigma \not= 0$, the coherent states do not exist. Let  $\psi_{\lambda , q}(N)$ be equal to $\lambda q^{-N}, q \leq 1 , \lambda >0$. When $\lambda = (q^{-1}-q)$, we get the function $\psi_{qosc}$ for $\sigma =0$ previously introduced and corresponding to the q-oscillator algebra (\ref{qosc}) with non-Fock representation and with coherent states.  Let us remark that this function also corresponds to $a a^\dagger = q^{-1}a^\dagger a$.

\nn In this case, the functional equations (\ref{mel1}):

\be
\hat{F}_{\lambda , q}(\rho+1)= \lambda q^{-\rho} \hat{F}_{\lambda , q}(\rho) 
\ee

\nn and
 
\be
x F_{\lambda , q}(x) = \lambda F_{\lambda , q}(qx)
\label{Fq}
\ee

\nn are equivalent for $F_{\lambda , q}(x)$ has to be fast decreasing at infinity and at the origin .

\nn Introducing $\ln \hat{F}_{\lambda , q}(\rho)$, we solve the resulting equation in terms of first and second Bernouilli polynomials so that we get the particular solution :

\be
\hat{F}^0_{\lambda , q}(\rho)= \exp\left( \rho \ln\lambda -\frac{1}{2}(\rho^2-\rho)\ln q \right)
\ee

\nn Taking the inverse Mellin transform, we obtain the following particular solution of (\ref{Fq}) \cite{kowalski} \cite{canada}:

\be
F^0_{\lambda , q}(x) = \exp\left(\frac{\ln^2(x\lambda ^{-1})}{2\ln (q)}-\frac{\ln(x\lambda ^{-1})}{2}\right)
\label{F0}
\ee

\nn The general solution of (\ref{Fq}) is given by :
\be
F_{\lambda , q}(x) = F^0_{\lambda , q}(x) h_{\lambda , q}(x)
\label{Fx}
\ee

\nn where $h_{\lambda , q}(x)$ is a function satisfying

\be
h_{\lambda , q}(x) = h_{\lambda , q}(qx)
\label{h}
\ee

\nn that is a periodic function of $\ln x/\ln q$ of period 1.
We verify directly that up to a constant the momentums $M(n)$ of $F^0_{\lambda , q}(x)$ are well equal to $\lambda^{-n} q^{-n(n+1)/2}$ for $n\in Z$, as wanted, and this is yet true for $F_{\lambda , q}(x)$ given in (\ref{Fx}) as implied by the definition of $h_{\lambda , q}(x)$.

\nn As the absolute value of a solution of equation (\ref{h}) is also a solution, the function $h_{\lambda , q}(x)$ can always be taken strictly positive.  Nevertheless, there exist an infinite number of $h_{\lambda , q}(x)$, in particular all the powers of any solution,  so that we get a lot of candidates to define the norm we are looking for. But mutatis mutandis, the situation is analogous to that initially studied by Bargmann and we can extend the largest part of its proof.

\nn Particularly we can prove that the set $\cal{S}$ of holomorphic functions under consideration is the set of functions $f(z)$ such that :

\be
\begin{array}{ll}
\mid f(z) \mid \leq C & \exp\left(-\frac{\ln^2(x\lambda ^{-1})}{4\ln q}-\frac{\ln(x\lambda ^{-1})}{4}\right) \times\\
&\sum _{-\infty}^{+\infty} \exp\left(\frac{\ln q}{2}(\frac{\ln (x\lambda^{-1})}{\ln q}+n+\frac{1}{2})^2 \right)
\end{array}
\ee

\nn The representation Hilbert space is the Hilbert space defined on $\cal{S}$ by

\be
(f,g)_F = \int F_{\lambda , q}(z\zz)  g(z) \overline{f(z)} dz d\zz
\label{fg}
\ee

\nn where $F_{\lambda , q}(z\zz)$ is a positive function of the form given in equation (\ref{Fx}). The various norms such obtained are proportional as verified on the basis elements.

\nn Furthermore, following always Bargmann's procedure, we can prove the closeness of the operator $z$ and $z^{-1}q^z\frac{d}{dz}$. Moreover  we have :

\be
\mid z f(z) \mid ^2 = \sum_{-\infty}^{\infty} \mid f_n \mid ^2  \psi (n+1)
\label{af}
\ee

\nn and

\be
\mid z^{-1}\psi ( z \frac {d}{dz}) f(z) \mid ^2 = \sum_{-\infty}^{\infty} \mid f_n \mid ^2  \psi (n)
\label{adagf}
\ee

\nn Using the two previous equations, we obtain :

\be
\mid z f(z)\mid ^2 = q \mid z^{-1} q^{z\frac{d}{dz}} f(z) \mid ^2
\ee

\nn  This proves that the operators $a$ and $a^\dagger$ have the same domain of definition and are  mutually adjoint  \cite{Bargmann}.

\nn We have proved the existence of  Bargmann representations for the deformed oscillator defined by (\ref{a1}) and (\ref{a3}) when $\psi_{q^{-1}-q , q} (N) = (q^{-1}-q) q^{-N}$. The resolution of identity is not unique and is given explicitly by (\ref{1}), (\ref{F0}) and (\ref{Fx}).

\subsection{Generalization of the previous example} 
 
\nn In this section, we point out some directions to extend the results of the previous section
when $\psi (x)$ is of the form :

\be
\psi (x) = \exp \left( \sum_0 ^{2p+1} a_n x^n \right), \quad a_{2p+1}>0
\ee

\nn A solution of Equation (\ref{mel1}) is :

\be
\hat{F}(\rho ) = \exp \left( \sum_0 ^{2p+1} \frac{a_n}{n+1}B_{n+1}(\rho)
  \right) 
\label{B}
\ee

\nn where $B_{n+1}(\rho)$ are the Bernouilli polynomials \cite{grad} defined by the difference equation:

\be
B_n(x+1)-B_n(x) =nx^{n-1}
\ee

The term of highest degree in (\ref{B}) is $a_{2p+1} \rho^{2p+2}/(2p+1)$. When $\rho$ is  pure imaginary $\rho = i \sigma , \sigma \in R$, this term is $a_{2p+1}(-1)^{p+1}\sigma^{2p+2}/(2p+1)$, therefore the inverse Mellin transform of $\hat{F}(\rho)$ exists only if $p$ is an even number. The function $F(x)$ thus obtained is always real, but not necessarily positive. Nevertheless, in specific cases, for example when the exponent in (\ref{B}) contains only the term of highest degree, $F(x)$ is actually strictly positive and the Bargmann procedure works as before.

\nn Since $ \psi(n+1)/\psi(n)$ grows indefinitely as $n \rightarrow \pm \infty$, Equations (\ref{af}) and (\ref{adagf}) imply that the domain of $z^{-1}\psi ( z \frac {d}{dz})$ is included in the domain of $z$ but cannot be identical. Nevertheless, the mutual adjointness can be proved on the basis, using equation (\ref{M}). 

\nn It is worthwhile to underline that this generalization gives an example where the Bargmann representation does not exist, namely when $p$ is odd.

\subsection {The ring case}

\nn Let us consider the simple case
 $\psi (x) = a + q^x , \quad q>1,\quad a> 0$, that is mainly involved in the study of the q-oscillator (\ref{qosc2}) with non-Fock representations. The domain of existence of the coherent states is $a\leq \mid z \mid^2$. The momentum $\hat{F}(n)$ reads :

\be
 \hat{F}(n)= \int_a ^{+\infty}F(x) x^{n-1}dx 
\ee

\nn We first prove that equations (\ref{M}) and (\ref{FF}) are not equivalent if $F(x) $ is positive on the whole positive axis. Indeed, let us start with a solution of (\ref{FF}), Equation (\ref{M}) reads :

\be
\int_{q^{-1}a}^{a} F(x) x^{n-1} dx =0
\ee

\nn that is obviously impossible if $F(x)$ is positive on $[q^{-1}a,a]$.
 Therefore in this case, the momentums deduced from the weight function solution of (\ref{FF}) are not the expected ones (solutions of (\ref{M})). 
Moreover, in \cite{canada}, we proved that the solution of (\ref{FF}) is identically zero.

\nn Let us look for a  solution of(\ref{mel1}) that cannot have poles, due to (\ref{nn}) :

\be
\hat{F}(\rho+1)=(q^\rho +a)\hat{F}(\rho)
\ee

\nn We have as a convenient particular solution the following entire function :

\be
\hat{F}(\rho)=  a^{\rho} \prod_{p \geq 0}(1+ a^{-1} q^{\rho-p-1})
\label{aa}
\ee

\nn but it is not a Mellin transform of a true function $F(x)$. Indeed if it has an inverse Mellin transform, it can be calculated on any parallel to the imaginary, for instance on $Re \rho =\ln a/\ln q$. On this axis, $\mid \hat{F}(iy)\mid \geq \prod_{p\geq 0}(1-q^{-p-1})$, so that (\ref{aa}) is not the Mellin transform of a true function. 

\nn Nevertheless, we can write (\ref{aa})  in the form \cite{exton}:

\be
\hat{F}(\rho) = a^{\rho}\left(1+\sum_{n\geq1} \frac{a^{-n}q^{n\rho}}{(q-1) \cdots (q^n-1)}\right)
\ee

\nn The series is absolutely convergent as $q>1$. It is easily verified that this expression can be seen as the Mellin transform of the following measure :

\be
F(x)= \sum_{n\geq 0} \frac{a^{-n}}{(q-1) \cdots (q^n-1)} \delta (\ln a+\ln q -\ln x)
\ee

\nn Therefore, in this case we obtain a Bargmann representation if we accept the weight function to be a true measure.

\nn The same is true when we consider $\psi(x) = 1/(q^x+a), q<1$. Equation (\ref{mel1}) reads :

\be
\hat{F}(\rho+1) = \frac{1}{q^\rho +a} \hat{F}(\rho)
\ee

The domain of existence of the coherent states is the disc of radius $1/a$. We then obtain :

\be
F(x)=\sum_{n\geq 0} \frac{q^{-n} a^{-n}}{(q^{-1}-1) \cdots (q^{-n}-1)}\delta (-\ln a+n\ln q -\ln x)
\ee 

\nn In this subsection, we gave examples where the Bargmann representations only exist if we admit that the scalar product be expressed bu means of true measures.

\section {$\psi$ vanishes } 

\nn In this section, we consider two cases where the spectrum of $N$ is the set $N^+$ of the positive integers and  where the coherent states are defined in the whole complex plane.

\subsection{q-oscillators}

 \nn The first example corresponds to (\ref{qosc})  with $\sigma =1$ and the second one to (\ref{qosc2}) with $\sigma = (q-1)^{-1}$.
\vspace{5mm}

\nn a)$\psi (x) = [x] \equiv (q^x -q^{-x})/(q-q^{-1})$
\vspace{5mm}

\nn A resolution of the identity was shown to be obtained with a q-integration \cite{gray}. The q-integration $\int_0^x d_q x$ is the inverse operator  of the q- derivative $D_q = \frac{1}{x} \frac{q^{x\partial _x }-q^{-x\partial _x}}{q-q^{-1}}$ that vanishes at the origin :

\be
\int_0^x d_q x = \frac{q-q^{-1}}{q^{x\partial_x}-q^{-x\partial_x}} x = (q^{-1}-q) \sum _{n\geq 0} q^{(2n+1)x\partial_x}x , \mbox{when} q<1
\ee

The q-exponential is defined by :

\be
Exp_q(x) = D_q Exp_q(x) = \frac{Exp_q(qx)-Exp_q(q^{-1}x)}{x(q-q^{-1})}
\label{expq}
\ee

\nn with the condition that it is equal to one when $x$ is zero. This
 function reads :

 \be
Exp_q(x)=\sum_{n\geq 0}\frac{x^n}{[n]!}
\ee

\nn and  vanishes on the negative axis \cite{gray}. Denoting by $-\zeta$ the first zero at the left of the origin
the resolution of identity then reads :

\be
\int \frac{d\theta}{2} \int_0^{\zeta^2} d_q \rho^2  Exp_q (-\rho^2)\mid \rho e^{-i\theta}>< \rho e^{i\theta}\mid 2 \rho d\rho =1
\label{11}
\ee

\nn Here we look for a Bargmann representation where the scalar product involves a true integral.

\nn First, it is easy to verify that in this case as in the following, if $F$ verifies (\ref{FF}), its Mellin transform is solution of (\ref{mel1}) and the moments are the expected ones. In both cases, we choose to define the weight function, not through its Mellin transform but directly as   solution of (\ref{FF}).

\nn The equation  (\ref{FF}) for this particular case reads  :

\be
x F(x) = \frac{q^{-x\partial _x}-q^{x\partial _x }}{q-q^{-1}}F(x)
\ee

\nn The obvious solution of this equation

\be
F(x)=Exp_q(-x)
\ee

\nn is not positive for all positive values of $x$  and the Bargmann representation as defined in section 4 does not exist. 

\nn Following the trick used to get (\ref{11}), we can try to limit the integration to the domain where $Exp_q(x)$ is positive. Let us see if 

\be
\hat{F}(n) = \int_0^{\zeta}Exp_q(-x) x^{n-1} dx
\ee

\nn could work. Equation (\ref{M}) gives :

\be 
\int_{q\zeta}^{\zeta} Exp_q(-q^{-1}x)x^{n-1} dx - \int_{q^{-1}\zeta}^{\zeta} Exp_q(-q x)x^{n-1} dx =0
\label{truc}
\ee

\nn The problem is symmetric under the change $q$ into $q^{-1}$. Let us choose $q > 1$, (\ref{truc}) takes the form :

\be 
 \int_{q^{-1}\zeta}^{\zeta}( Exp_q(-x)q^n+ Exp_q(-q x)) x^{n-1} dx =0
\label{truc1}
\ee 

\nn The integrand of (\ref{truc1}) reads $( Exp_q(-x)(q^n -x(q-q^{-1}))+ Exp_q(-q^{-1} x)) x^{n-1}$ and is positive for $n$ enough large; this leads to $Exp_q(-x)(q^n -x(q-q^{-1})+ Exp_q(-q^{-1} x)=0$, that is impossible. 

\nn Therefore in order to obtain a Bargmann representation, we must look for another solution of (\ref{FF}) that will be positive. As already noticed, the problem being symetric under the change $q$ into $q^{-1}$, we assume $q>1$. Let us start with (\ref{mel1})that reads :

\be
\hat{F}(\rho+1)=\frac{q^\rho-q^{-\rho}}{q-q^{-1}} = \frac{q^\rho}{q-q^{-1}}(1-q^{-2\rho})\hat{F}(\rho)
\ee

\nn Let us write $\hat{F}$ on the form :

\be
\hat{F}(\rho)= \phi q^{\frac{\rho}{2}(\rho -1)}(q-q^{-1})^{-\rho} \hat{f}(\rho)
\label{fr}
\ee

\nn The function $\hat{f}(\rho)$ must verify

\be
\hat{f}(\rho+1)=(1-q^{-2\rho})\hat{f}(\rho)
\ee

\nn and is given by :

\be
\hat{f}(\rho)= \sum_{n\geq 0}\frac{q^{-2n\rho}}{(1-q^{-2}) \cdots (1-q^{-2n})}
\label{fr1}
\ee

\nn The condition $\hat{F}(0)=1$, furnishes the normalization factor :

\be
\phi = (q-q^{-1})\left( 1+\sum_{n>0}\frac{q^{-2n}}{(1-q^{-2}) \cdots (1-q^{-2n})}\right)^{-1}
\label{fr2}
\ee

\nn Putting (\ref{fr1}) and (\ref{fr2}) in (\ref{fr}), we obtain $\hat{F}(\rho)$ 
, and then we can calculate its inverse Mellin transform :

\be
F(x)= \frac{\exp(-\frac{1}{2\ln q}(\ln x +\ln(q-q^{-1})+\frac{1}{2}\ln q)^2)}{ 1+\sum_{n>0}\frac{q^{-2n}}{(1-q^{-2}) \cdots (1-q^{-2n})}}\sum_{n\geq 0}\frac{q^{-n(2n+1)}((q-q^{-1})x)^{-2n}}{(1-q^{-2}) \cdots (1-q^{-2n})}
\ee

\nn This function being positive, we have obtained a Bargmann representation where the scalar product is written with a true integral. Let us stress that $F(-x)$ is solution of (\ref{expq}) and is thus a possible candidate to write the resolution of identity with a q-integration and a positive function on the whole positive axis.

\nn The same is true for in the next example where two resolutions of the identity coexist.
\vspace{5mm}

\nn b) $\psi (x) =(x) \equiv (q^x-1)/(q-1)$, with $ q > 1$

\vspace{5mm}
\nn First  we show that the resolution of the identity can be obtained with a q-integral as in \cite{gray}. 

\nn The q-integration \cite{jackson}, \cite{exton}, \cite{macanally} is defined to be the inverse of the q-derivative $D_q \equiv \frac{1}{x} \frac{q^{x\partial _x }-1}{q-1}$  :

\be
\int_0^x d_q x  \equiv \frac{q-1}{q^{x\partial_x}-1} x = (q-1) \sum _{n \geq 0 }q^{-(n+1)x\partial_x}x
\ee

\nn The q-exponential, solution of the equation :
\be
 Exp_q(x) = D_q Exp_q(x)= \frac{Exp_q(qx)-Exp_q(x)}{x(q-1)} 
\ee

\nn is given by :

 \be
Exp_q(x) = \prod_{p\ge 0} (1+x(1-q^{-1})q^{-p})= \sum_{n\geq 0}\frac{x^n}{(n)!}
\ee

\nn and vanishes for $x = - q^p (1-q^{-1})^{-1}$. The nearest zero on the left of the origin is $-\zeta = -(1-q^{-1})^{-1}$. Therefore 
the resolution of the identity takes the same form as in (\ref{11}) with the new expressions for $\int_0^x d_q x$, $Exp_q$ and $\zeta$.

\nn Let us now look for a Bargmann representation as defined in section 4. We see that  the equation (\ref{FF}) can be written :

\be
F(q^{-1}x)= ( x(q-1)+1) F(x)
\ee

\nn We  easily prove that the weight function is given by 

\be
F(x)= \frac{1}{Exp_q(qx)}
\ee
 
\nn It is a positive function when $x >0$ and  its Mellin transform fulfills :

\be
 \hat{F}(\rho+1)=\frac{q^\rho -1}{q-1}\hat{F}(\rho)
\ee

\nn This ensures that the momentum $\hat{F}(n)$ are the expected one (\ref{moments}). Thus, in this case, coexist two resolutions of the identity ,one involving a true integral and a weight function $F(x) = (Exp_q(qx))^{-1}$ and one with a q-integral, the weight function being $Exp_q(-x)$.

\subsection{ $\psi (x) = x^n , n> 0$}

\nn The Mellin transform of the weight function is solution of the equation deduced from (\ref{mel1}) :

\be
\hat{F}(\rho +1) = \rho ^n \hat{F}(\rho)
\ee

\nn and can be expressed with the gamma-function $\Gamma (z) = \int_0^\infty e^{-t} t^{z-1} dt$ :

\be
\hat{F}(\rho)= \left( \Gamma (\rho) \right )^n
\ee

\nn When $n$ is an integer, the inverse Mellin transform gives $F(x)$

\be
F(x) =\int_0^\infty \cdots \int_0^\infty e^{-(t_1+\cdots + t_n)} dt_1 \cdots dt_n \delta (x-t_1 \times \cdots \times t_n)
\ee

\nn On this expression, we see that $F(x) $ is a positive function so that the Bargmann representation exists. In the case $n=1$, we recover the usual harmonic oscillator where $F(x)= e^{-x}$.

\section{Conclusion}

\nn We have studied the possibility of Bargmann representations for any deformed oscillator algebra characterized by a function $\psi$. We gave the conditions to be verified by this function for admitting representations with coherent states. We get the unique functional equation to be satisfied by the Mellin transform of the weight function defining the scalar product. We were able to get definite and positive answer in many cases including in particular some types of q-oscillators. Although we cannot succeed to obtain a general characterization of the function $\psi$ leading to Bargmann representations, we underline two points : 

\nn - We exhibit cases where the Bargmann representations do not exist even when coherent states do ( subsection (5.2));

\nn - The analysis of subsection (5.3) shows that the scope of our study have to be extended up to include true measures for writing the scalar product.

\nn Finally let us remark that we have obtained  scalar products for the Bargmann representations of the usual q-oscillators, involving  true integrals instead of  q-integrations as previously proposed in literature.

\end{document}